\documentclass[12pt,leqno]{article}
\usepackage{graphicx}
\title{Multitime maximum principle approach of
minimal submanifolds and harmonic maps}
\author{Constantin Udri\c ste}

\date{}
\begin{document}
\maketitle
\textwidth=14cm
\textheight=21cm

\newcommand{\di}{\displaystyle}
\newcommand{\ov}{\over}
\newcommand{\noa}{\noalign{\medskip}}
\newcommand{\al}{\alpha}
\newcommand{\be}{\beta}
\newcommand{\om}{\omega}
\newcommand{\bt}{\bar\tau}
\newcommand{\br}{\hbox{\bf R}}
\newcommand{\bd}{\hbox{\bf D}}
\newcommand{\bq}{\hbox{\bf Q}}
\newcommand{\bn}{\hbox{\bf N}}
\newcommand{\bz}{\hbox{\bf Z}}
\newcommand{\bc}{\hbox{\bf C}}
\newcommand{\sgn}{\hbox{sgn}}
\newcommand{\ome}{\Omega}
\newcommand{\ba}{\bar a}
\newcommand{\ti}{\times}
\newcommand{\ga}{\gamma}
\newcommand{\ty}{\infty}
\newcommand{\de}{\delta}
\newcommand{\te}{\theta}
\newcommand{\pa}{\partial}
\newcommand{\va}{\varphi}
\newcommand{\ld}{\ldots}
\newcommand{\qu}{\quad}
\newcommand{\la}{\lambda}
\newcommand{\fo}{\forall}
\newcommand{\ep}{\varepsilon}
\newcommand{\pp}{\prime}
\newcommand{\su}{\subset}
\newcommand{\si}{\sigma}
\newcommand{\gaa}{\Gamma}
\newcommand{\ri}{\Rightarrow}
\newcommand{\rii}{\Leftrightarrow}
\newcommand{\mm}{\medskip}
\newcommand{\na}{\nabla}

\abstract{Some optimization problems coming from the Differential Geometry, as for example,
the minimal submanifolds problem and the harmonic maps problem are solved here via interior solutions of
appropriate multitime optimal control problems.

Section 1 underlines some science domains where appear multitime
optimal control problems. Section 2 (Section 3) recalls the
multitime maximum principle for optimal control problems with
multiple (curvilinear) integral cost functionals and $m$-flow type
constraint evolution. Section 4 shows that there exists a multitime
maximum principle approach of multitime variational calculus.
Section 5 (Section 6) proves that the minimal submanifolds (harmonic
maps) are optimal solutions of multitime evolution PDEs in an
appropriate multitime optimal control problem. Section 7 uses the
multitime maximum principle to show that of all solids having a
given surface area, the sphere is the one having the greatest
volume. Section 8 studies the minimal area of a multitime linear
flow as optimal control problem. Section 9 contains commentaries.}



\section{Origin of multitime optimal control problems}

The original results in this paper show that the minimal submanifolds (see [1], [3], [4], [22], [31]) and the harmonic maps (see [2])
are solutions of multitime optimal control problems (see [5]-[30]). In this way we change the traditional geometrical viewpoint,
looking at a minimal submanifold or at a harmonic map as solution in a multitime optimal control system via the multitime maximum principle.
Problems of this type are among the most challenging in Differential Geometry and Control Theory while being among the most important for various applications.

Near Differential Geometry, there are many domains of science
containing multitime optimal control problems: Material Strength (the description of torsion of prismatic bars in
the elastic case as well as in the elastic-plastic case),
Fluid Mechanics (the motion of  fluid substances, Navier-Stokes PDE written as first order PDE),
Magnetohydrodynamics (Maxwell-Vlasov PDE, Navier-Stokes PDE) etc.
Of course, to describe some $m$-dimensional objects as optimal evolution maps,
a deeply understanding of the meaning of evolution is necessary.
Our results confirm that the basic ideas of Lev S. Pontryaguin, Lawrence C. Evans and Jacques-Louis Lions
are still surviving for optimal control problems governed by normal first order PDEs.

\section{Optimal Control Problem with \\Multiple Integral Cost Functional}

Let us analyze a multitime optimal control problem based on a multiple integral cost functional and $m$-flow type $PDE$ constraints [14]-[15],
[18], [19], [28], [30]:
$$
\di{\max_{u(\cdot)}}\,\,I(u(\cdot))=\int_{\Omega_{0t_0}}L(t,x(t),u(t))dt
$$
$$
\hbox{subject to}
$$
$$
\frac{\partial x^i}{\partial t^\alpha}(t) = X^i_\alpha (t,x(t),u(t)), i=1,...,n; \alpha = 1,..., m,
$$
$$
u(t)\in {\cal U}, \,\, t\in {\Omega_{0t_0}};\,\, x(0)=x_0,\,\, x(t_0)=x_{t_0}.
$$
Mathematical data: $t = (t^\alpha)=(t^1,...,t^m)\in R^m_+$ is the {\it multitime} (multi-parameter of evolution); $dt=dt^1\wedge...\wedge dt^m$
is the {\it volume element} in $R^m_+$; $\Omega_{0t_0}$ is the parallelepiped fixed by the diagonal opposite points
$0=(0,...,0)$ and $t_0=(t^1_0,...,t^m_0)$ which is equivalent to the {\it closed interval} $0\leq t\leq t_0$ via the {\it product order} on $R^m_+$;
$x: \Omega_{0t_0}\to R^n,\, x(t)=(x^i(t))$ is a $C^2$ {\it state vector};
$u: \Omega_{0t_0}\to U\subset R^k,\,u(t)=(u^a(t)),\,\, a=1,...,k$ is a $C^1$ {\it control vector};
the {\it running cost} $L(t, x(t),u(t))$ is a $C^1$ {\it nonautonomous Lagrangian}; $X_\al (t,x(t),u(t))=(X^i_\alpha (t,x(t),u(t)))$ are $C^1$
vector fields satisfying the complete integrability conditions ($m$-flow type problem), i.e., $D_\be X_\al = D_\al X_\be$
($D_\al$ is the total derivative operator) or
$$\left({\pa X_\al \ov \pa u^a} \de^\ga_\be - {\pa X_\be \ov \pa
u^a} \de^\ga_\al\right) {\pa u^a\ov \pa t^\ga} = [X_\al, X_\be] +\di\frac{\pa X_\be}{\pa t^\al}-\di\frac{\pa X_\al}{\pa t^\be},$$
where $[X_\al, X_\be]$ means the {\it bracket} of vector fields.
The complete integrability hypothesis constrains the set of all admissible controls (satisfying the complete integrability conditions)
$$
{\cal U} =\left \{\left. u: R^m_+ \to U\, \right| D_\be X_\al = D_\al X_\be\right\}
$$
and the admissible states.

To formulate the {\it multitime maximum principle} we need the {\it control Hamiltonian}
$$H(t,x(t),u(t),p(t))=L(t,x(t),u(t))+p_i^\al(t)X^i_\al(t,x(t),u(t)).$$

{\bf Theorem 1.1 (multitime maximum principle; necessary conditions)}.
{\it Suppose that the previous optimal control problem,
with $X, X^i_\al$ of class $C^1$, has an interior solution $\hat u(t)\in {\cal U}$ which
determines the $m$-sheet of state variable $x(t)$. Then there exists a $C^1$ costate
$p(t)=(p^\al_i(t))$ defined over  $\Omega_{0t_0}$ such that the relations
$$\di\frac{\pa p^\al_j}{\pa t^\al}(t) = - \di\frac{\pa H}{\pa x^j}(t,x(t),\hat u(t),p(t)), \,\, \forall t\in {\Omega_{0t_0}},$$
$${\delta_{\al\be}p^\al_j(t)n^\be(t) |}_{\pa \Omega_{0t_0}} = 0\,\hbox{(orthogonality or tangency)},$$
$$\di\frac{\pa x^j}{\pa t^\al}(t)=\di\frac{\pa H}{\pa p^\al_j}(t,x(t),\hat u(t),p(t)), \,\, \forall t\in {\Omega_{0t_0}},\,\,x(0)=x_0$$
and
$$ H_{u^a}(t,x(t),\hat u(t),p(t))=0, \,\, \forall t\in {\Omega_{0t_0}}$$
hold}.

{\bf Remark 1.2} If the optimal control $\hat u(t)\in {\cal U}$ is not an interior point,
then instead of critical point condition we have
$$H((t,x(t),\hat u(t),p(t)) = \max_{u}H(t,x(t),u,p(t)).$$

{\bf Theorem 1.2 Multitime Global Maximality}. {\it The foregoing problem has a global solution $(\hat x(\cdot),\hat u(\cdot))$ if and only if 
the multiple integral functional $\int_{\Omega_{0t_0}}H(t,x(t),u(t),p(t))dt$ is incave with respect to $u$.}

\section{Optimal Control Problem with \\Curvilinear Integral Cost Functional}

A multitime optimal control problem whose cost functional is the sum between a path independent curvilinear integral (mechanical work, circulation)
and a function of the final event, and whose evolution PDE is an $m$-flow, has the form [14], [15]-[17], [19], [23], [24], [30]
$$
\begin{array}{cc}
\hspace{0.cm}\di{\max_{u(\cdot)}}&\hspace{-1.3cm}J(u(\cdot))= \di\int_{\Gamma_{0t_0}}L_\al(t,x(t),u(t))dt^\al + g(x(t_0))\\ \noa
\hbox{subject to} \,\,&\di{\pa x^i \ov \pa t^\al}(t) = X^i_\al (t,x(t),u(t)),\,\, i=1,...,n,\,\, \al = 1,..., m, \\ \noa
&\hspace{-1.3cm}u(t)\in {\cal U}, \,\,t\in {\Omega_{0t_0}}, \,\,x(0)=x_0,\,\, x(t_0)=x_{t_0}.
\end{array}
$$

This problem requires the following data: the {\it multitime} (multiparameter of evolution) $t = (t^\al)\in R^m_+$; an
arbitrary $C^1$ curve $\Gamma_{0t_0}$ joining the diagonal opposite points
$0=(0,...,0)$ and $t_0=(t^1_0,...,t^m_0)$ in the parallelepiped $\Omega_{0t_0}=[0,t_0]$
(multitime interval) in $R^m_+$ endowed with the {\it product order};
a $C^2$ {\it state vector} $x:\Omega_{0t_0}\to R^n,\,x(t)=(x^i(t))$; a $C^1$ {\it control vector}
$u: \Omega_{0t_0}\to U\subset R^k,\, u(t)=(u^a(t)),\,\, a=1,...,k$;
a {\it running cost} $L_\al(t, x(t),u(t))dt^\al$ as a nonautonomous closed (completely integrable) {\it Lagrangian 1-form}, i.e.,
it satisfies $D_\be L_\al = D_\al L_\be$ ($D_\al$ is the total derivative operator) or
$$\left({\pa L_\al \ov \pa u^a} \de^\ga_\be - {\pa L_\be \ov \pa
u^a} \de^\ga_\al\right) {\pa u^a\ov \pa t^\ga} = X^i_\al\di\frac{\pa L_\be}{\pa x^i} - X^i_\be\di\frac{\pa L_\al}{\pa x^i}
+\di\frac{\pa L_\be}{\pa t^\al}-\di\frac{\pa L_\al}{\pa t^\be};$$
the {\it terminal cost functional} $g(x(t_0))$; the $C^1$ vector fields $X_\al=(X^i_\al)$ satisfying
the complete integrability conditions ($m$-flow type problem), i.e., $D_\be X_\al = D_\al X_\be$ or
$$\left({\pa X_\al \ov \pa u^a} \de^\ga_\be - {\pa X_\be \ov \pa
u^a} \de^\ga_\al\right) {\pa u^a\ov \pa t^\ga} = [X_\al, X_\be] +\di\frac{\pa X_\be}{\pa t^\al}-\di\frac{\pa X_\al}{\pa t^\be},$$
where $[X_\al, X_\be]$ means the {\it bracket} of vector fields.
Some of the previous hypothesis select the set of all admissible controls (satisfying the complete integrability conditions,
eventually, a.e.)
$$
{\cal U} =\left \{\left. u: R^m_+ \to U\, \right| D_\be L_\al = D_\al L_\be,\,\,D_\be X_\al = D_\al X_\be,\,\,\hbox{a.e.}\right\}.
$$
The set ${\cal U}$ does not contain always the constant controls, but it contain sure controls which are continuous at the right
(in the sense of product order).

The previous PDE evolution system is equivalent to the path-independent curvilinear integral equation
$$x(t)=x(0) +\di\int_{\ga_0t} X_\al(s,x(s),u(s))ds^\al,$$
where $\ga_{0t}$ is an arbitrary piecewise $C^1$ curve joining the opposite diagonal points $0$ and $t$
of the parallelepiped $\Omega_{0t}=[0,t]\subset \Omega_{0t_0}=[0,t_0]$.

In the multitime optimal control problems with path independent integrals, it is enough to use increasing curves.

{\bf Definition 2.1} A piecewise $C^1$ curve $\ga_{0t_0}: s^\al = s^\al (\tau),\,\,\tau \in [\tau_0,\tau_1]$, $s(\tau_0)=0, s(\tau_1)=t_0$
is called {\it increasing} if the tangent vector $\dot s =(\dot s^\al)$ satisfies $\dot s^\al \geq 0$, with $||\dot s||=0$ only at
isolated points.

If we use the {\it control Hamiltonian 1-form}
$$H_\al(t,x(t),u(t),p(t))=L_\al(t,x(t),u(t))+p_i(t)X^i_\al(t,x(t),u(t)),$$
we can formulate the {\it multitime maximum principle}

{\bf Theorem 2.1 Multitime Maximum Principle}. {\it Suppose that the previous problem,
with $L_\al, X^i_\al$ of class $C^1$, has an interior solution $\hat u(t)\in {\cal U}$ which determines the
$m$-sheet of state variable $x(t)$. Then, there exists a $C^1$ vector costate
$p(t)=(p_i(t))$ defined over  $\Omega_{0t_0}$ such that the following relations hold:
$$\di\frac{\pa p_j}{\pa t^\al}(t) = - \di\frac{\pa H_\al}{\pa x^j}(t,x(t),\hat u(t),p(t)),\,\,\forall t\in {\Omega_{0t_0}},\,\, p_j(t_0) = 0,$$
$$\di\frac{\pa x^j}{\pa t^\al}(t)=\di\frac{\pa H_\al}{\pa p_j}(t,x(t),\hat u(t),p(t)),\,\,\forall t\in {\Omega_{0t_0}}, \,\,x(0)=x_0,$$
$$ \di\frac{\pa H_\al}{\pa u^a}(t,x(t),\hat u(t),p(t))=0, \,\, \forall t\in {\Omega_{0t_0}}.$$}

{\bf Remark 2.1} If the optimal control $\hat u(t)\in {\cal U}$ is not an interior point,
then instead of critical point condition we have
$$H_\al((t,x(t),\hat u(t),p(t)) = \max_{u}H_\al(t,x(t),u,p(t)).$$

{\bf Theorem 2.2 Multitime Global Maximality}. {\it The foregoing problem has a global solution $(\hat x(\cdot),\hat u(\cdot))$ if and only if 
the curvilinear integral functional $\di\int_{\Gamma_{0t_0}}H_\al(t,x(t),u(t),p(t))dt^\al$ is incave with respect to $u$.}

\section{Multitime maximum principle \\approach of variational calculus}

In fact we show that the {\it multitime maximum principle} motivates the {\it multitime Euler-Lagrange or Hamilton PDEs}.

\subsection{Case of multiple integral action}

Suppose that the evolution system is reduced to a completely integrable system
$$\di\frac{\pa x^i}{\pa t^\al}(t) = u^i_\al (t),\,\, x(0)=x_0,\,\, t \in \Omega_{0t_0}\subset R^m_{+},\eqno (PDE)$$
and the functional is a multiple integral
$$
I(u(\cdot))=\int_{\Omega_{0t_0}}X(t,x(t),u(t))dt, \eqno(I)
$$
where $\Omega_{0t_0}$ is the parallelepiped fixed by the diagonal opposite points
$0=(0,...,0)$ and $t_0=(t^1_0,...,t^m_0)$, the {\it running cost}
$X (t,x(t),u(t))dt$ is a {\it Lagrangian $m$-form}.

The associated basic control problem leads necessarily to the
multitime maximum principle. Therefore, to solve it we need the
control Hamiltonian
$$H(t,x,p_0,p,u) =  X(t,x,u) + p_i^\al u^i_\al$$
and the adjoint PDEs
$$\di\frac{\pa p_i^\al}{\pa t^\al}(t) = - \di\frac{\pa X}{\pa x^i}(t,x(t),u(t)). \eqno (ADJ)$$
Suppose the simplified multitime maximum principle is applicable
$$\di\frac{\pa H}{\pa u^i_\ga} = \di\frac{\pa X}{\pa u^i_\ga} + p_i^\ga = 0,\,\,
p_i^\ga = -\di\frac{\pa X}{\pa u^i_\ga},\,\, u^i_\ga = x^i_\ga.$$
Suppose the function $X$ is dependent on $x$ (a strong condition!).
We eliminate $p^\al_i$ using the adjoint PDE. It follows the multitime Euler-Lagrange PDEs
$$\di\frac{\pa X}{\pa x^i} - D_\al\left(\di\frac{\pa X}{\pa x^i_\al}\right)=0.$$

\subsection{Case of path independent curvilinear integral action}

First, suppose that the evolution system is reduced to a completely integrable system
$$\di\frac{\pa x^i}{\pa t^\al}(t) = u^i_\al (t),\,\, x(0)=x_0,\,\, t \in \Omega_{0t_0}\subset R^m_{+},\eqno (PDE)$$
and the functional is a path independent curvilinear integral
$$
J(u(\cdot))= \di\int_{\Gamma_{0t_0}} L_\be (t,x(t),u(t))dt^\be, \,\,u =(u^i_\al), \eqno (J)
$$
where $\Gamma_{0t_0}$ is an arbitrary piecewise $C^1$ curve joining the points $0$ and $t_0$, the {\it running cost}
$\omega = L_\be (t,x(t),u(t))dt^\be$ is a closed (completely integrable) {\it Lagrangian $1$-form}.

The associated basic control problem leads necessarily to the
multitime maximum principle. Therefore, to solve it we need the
control Hamiltonian $1$-form
$$H_\be(t,x,p_0,p,u) =  L_\be(t,x,u) + p_i u^i_\be$$
and the adjoint PDEs
$$\di\frac{\pa p_i}{\pa t^\be}(t) = - \di\frac{\pa L_\be}{\pa x^i}(t,x(t),u(t)). \eqno (ADJ)$$
Suppose the simplified multitime maximum principle is applicable, i.e.,
$$\di\frac{\pa H_\be}{\pa u^i_\ga} = \di\frac{\pa L_\be}{\pa u^i_\ga} + p_i \de ^\ga_\be = 0,\,\,
p_i \de ^\ga_\be = -\di\frac{\pa L_\be}{\pa u^i_\ga},\,\, u^i_\ga = x^i_\ga.$$
We accept that the functions $L_\be$ are dependent on $x$ (a strong condition!). Then $(ADJ)$ shows that
$$p_i(t) = p_i(0) - \di \int_{\Gamma_{0t}}{\pa L_\be \ov \pa x^i}(s,x(s),u(s))ds^\be,$$
where $\Gamma_{0t}$ is an arbitrary piecewise $C^1$ curve joining
the points $0, t\in \Omega_{0t_0}$. From the foregoing last three relations, it
follows
$$-\di\frac{\pa L_\be}{\pa x^i_\ga}(t,x(t),u(t))=\de^\ga_\be p_i(0) - \de^\ga_\be \di \int_{\Gamma_{0t}}\di\frac{\pa L_\la}{\pa x^i}(s,x(s),u(s))ds^\la.$$
If $X^0_\be$ are functions of class $C^2$, then applying the divergence operator
$D_\ga = \di\frac{\pa}{\pa t^\ga}$, we find {\it the multitime Euler-Lagrange $PDEs$}
$$\di\frac{\pa L_\be}{\pa x^i} - D_\ga\di\frac{\pa L_\be}{\pa x^i_\ga} = 0.$$

\section{Minimal submanifolds as \\optimal evolutions}

The minimal submanifolds are characterized by zero mean curvature. These
become an area of intense mathematical and scientific study over the past 15 years,
specifically in the areas of molecular engineering and materials sciences
due to their anticipated nanotechnology applications. The most extensive meeting ever held on the subject, in its 250-year history,
was organized in 2001 at Clay Mathematics Institute. In spite of all these efforts, the traditional thinking about minimality was not change.

Recently, I gave a new approach to the theory of minimal surfaces [22], [29] changing the traditional approach into a new one based on
solutions in a two-time optimal control system via the multitime maximum principle. More,
our recent studies [5]-[30] show the efficiency  of approaching some classical problems using techniques of multitime optimal control
or multitime modeling.

Let $\Omega_{0\tau}$ be an interval fixed by the diagonal opposite points $0, \tau \in R^m_+$.
On an $n$-dimensional Riemannian manifold $(M,g_{ij})$ we introduce the multitime controlled dynamics
$$
{\pa x^i \ov \pa t^\al}(t) = u^i_\al(t),\,\, i=1,...,n; \al=1,...,m; \eqno (PDE)
$$
$$
t=(t^1,...,t^m)\in \Omega_{0\tau},\,\,x^i(0)=x^i_0, \,\,x^i(\tau)=x^i_1,
$$
where $u= (u_\al) = (u^i_\al):\Omega_{0\tau}\to R^{mn}$ represents two open-loop $C^1$ control vectors,
linearly independent, eventually fixed on the boundary $\pa \Omega_{0\tau}$.
The complete integrability conditions of the (PDE) system restrict the set of controls to
$$
{\cal U} =\left \{\left. u=(u_\al) = (u^i_\al)\, \right| \di\frac{\pa u^i_\al}{\pa t^\be}(t)= \di\frac{\pa u^i_\be}{\pa t^\al}(t)\right\}.
$$
A $C^2$ solution of $(PDE)$ system is a submanifold (m-sheet) $\sigma: x^i=x^i(t^1,...,t^m)$.

The metric $g_{\al\be}(t) = g_{ij}(x(t))u^i_\al(t) u^j_\be(t)$ determines the volume $$\int_{\Omega_{0\tau}}\sqrt{\det(g_{\al\be})}\,\,dt$$
of the $m$-sheet $x(t),\,\,t\in \Omega_{0\tau}$ and this defines the functional
$$
V(u(\cdot)) = - \int_{\Omega_{0\tau}}\sqrt{\det(g_{\al\be})}\,\,dt. \eqno (V)
$$

{\bf Multitime optimal control problem of minimal submanifolds:}
$$
\di{\max_{u(\cdot)}}\,\,V(u(\cdot))
$$
$$
\hbox{subject to}
$$
$$
{\pa x^i \ov \pa t^\al}(t) = u^i_\al(t),\,\, i=1,...,n; \,\, \alpha = 1,..., m,
$$
$$
u(t)\in {\cal U}, \,\, t\in {\Omega_{0\tau}};\,\, x(0)=x_0,\,\, x(\tau)=x_1.
$$

To solve the multitime optimal control problem of minimal submanifolds, we apply the multitime maximum
principle. In general notations, we have
$$x = (x^i),\,\,u=(u_\al) = (u^i_\al),\,\,p=(p^\al) =(p^\al_i), \,\,i=1,...,n; \,\, \alpha = 1,..., m,$$
$$X_\al(x(t),u(t)) = u_\al(t),\,\,L(x(t),u(t)) = - \sqrt{\det(g_{\al\be})}$$
and the control Hamiltonian is
$$
H(x,u,p) = p^\al_iX^i_\al(x,u) + p_0 L(x,u).
$$
Taking $p_0=1$, the adjoint dynamics says
$$
{\pa p^\al_i \ov \pa t^\al} = -{\pa H\ov \pa x^i}.\eqno (ADJ)
$$
On the other hand, we have to maximize $H(x,u,p)$ with respect to the control $u$, hence
$$\di\frac{\pa H}{\pa u^i_\al} = \di\frac{\pa L}{\pa u^i_\al} + p_i^\al = 0.$$
This necessary condition (critical point) is also sufficient since the foregoing Hamiltonian $H(x,u,p)$ is an incave function with respect to $u$.
Indeed, the Lagrangian $L(x,u)$ is an incave function with respect to $u$, because it has any critical point, and the Hamiltonian $H$
is obtained from $L$ adding a linear term.

Having in mind that $x^i_\al = u^i_\al$, we eliminate $p^\al_i$
using the adjoint PDE. It follows the multitime Euler-Lagrange PDEs
$$\di\frac{\pa L}{\pa x^i} - D_\al\left(\di\frac{\pa L}{\pa x^i_\al}\right)=0$$
or
$$\di\frac{\pa L}{\pa x^i} - D_\al \left(\di\frac{\pa L}{\pa g_{\be\ga}}\,\,\di\frac{\pa g_{\be\ga}}{\pa x^i_\al}\right)=0.$$
Summarizing, we obtain

{\bf Theorem 5.1}. {\it A $C^2$ solution of the previous optimal control problem is a solution of the boundary problem
$$\di\frac{\pa L}{\pa x^i} - D_\al \left(g\,\, g^{\be\ga}\,\,\di\frac{\pa g_{\be\ga}}{\pa x^i_\al}\right)=0,\,\, x(0)=x_0,\,\, x(t_0)=x_{t_0}$$
i.e., it is a minimal submanifold.}

The familiar relativistic free particle and the dual string, which provides a dynamical model of hadrons, are described by the PDEs in Theorem 5.1.

\section{Harmonic maps as \\optimal evolutions}

Traditionally, the Harmonic maps are solutions to a natural geometric variational problem
motivated by some fundamental ideas from differential geometry,
in particular geodesics, minimal surfaces, and harmonic functions. Harmonic maps are also
closely related to nonlinear partial differential equations,
holomorphic maps in several complex variables, the theory of stochastic processes,
the nonlinear field theory in theoretical physics, and the theory of liquid crystals in materials science.

Our recent papers change the traditional geometrical viewpoint, looking at a harmonic map
as solution in a multitime optimal control system via the multitime maximum principle.

Let $\Omega_{0\tau}$ be a interval fixed by the diagonal opposite points $0, \tau \in R^m_+$.
Let $h_{\al\be}$ be a Riemannian metric on $R^m_+$.
On the Riemannian manifold $(M, g_{ij})$ we introduce the multitime controlled dynamics
$$
{\pa x^i \ov \pa t^\al}(t) = u^i_\al(t),\,\, i=1,...,n; \al=1,...,m; \eqno (PDE)
$$
$$
t=(t^1,...,t^m)\in \Omega_{0\tau},\,\,x^i(0)=x^i_0, \,\,x^i(\tau)=x^i_1,
$$
where $u= (u_\al) = (u^i_\al):\Omega_{0\tau}\to R^{mn}$ represents two open-loop $C^1$ control vectors,
linearly independent, eventually fixed on the boundary $\pa \Omega_{0\tau}$.
The complete integrability conditions of the (PDE) system, restrict the set of controls to
$$
{\cal U} =\left \{\left. u=(u_\al) = (u^i_\al)\, \right| \di\frac{\pa u^i_\al}{\pa t^\be}(t)= \di\frac{\pa u^i_\be}{\pa t^\al}(t)\right\}.
$$
A $C^2$ solution of $(PDE)$ system is a map (m-sheet) $\sigma: x^i=x^i(t^1,...,t^m)$.

We introduce the energy density $\frac{1}{2}h^{\al\be}(t)g_{ij}(x(t))u^i_\al(t) u^j_\be(t)$ of the $m$-sheet $x(t),\,\,t\in \Omega_{0\tau}$
and the energy functional (elastic deformation energy)
$$
E(u(\cdot)) = - \frac{1}{2}\int_{\Omega_{0\tau}}h^{\al\be}(t)g_{ij}(x(t))u^i_\al(t) u^j_\be(t)\,\,dt. \eqno (E)
$$
The energy density is defined by the trace of the induced metric $g_{\al\be}=g_{ij}u^i_\al u^j_\be$ with respect to the metric $h^{\al\be}$.

{\bf Multitime optimal control problem of harmonic maps:}
$$
\di{\max_{u(\cdot)}}\,\,E(u(\cdot))
$$
$$
\hbox{subject to}
$$
$$
{\pa x^i \ov \pa t^\al}(t) = u^i_\al(t),\,\, i=1,...,n; \,\, \alpha = 1,..., m,
$$
$$
u(t)\in {\cal U}, \,\, t\in {\Omega_{0\tau}};\,\, x(0)=x_0,\,\, x(\tau)=x_1.
$$

To solve the previous problem we apply the multitime maximum
principle. In general notations, we have
$$x = (x^i),\,\,u=(u_\al) = (u^i_\al),\,\,p=(p^\al) =(p^\al_i), i=1,...,n; \,\, \alpha = 1,..., m,$$
$$X_\al(x(t),u(t)) = u_\al(t),\,\,L(x(t),u(t)) = - \frac{1}{2} h^{\al\be}(t)g_{ij}(x(t))u^i_\al(t) u^j_\be(t)$$
and the control Hamiltonian is
$$
H(x,u,p) = p^\al_iX^i_\al(x,u) + p_0 L(x,u).
$$
Taking $p_0=1$, the adjoint dynamics says
$$
{\pa p^\al_i \ov \pa t^\al} = -{\pa H\ov \pa x^i}.\eqno (ADJ)
$$
On the other hand, we have to maximize $H(x,u,p)$ with respect to the control $u$, hence
$$\di\frac{\pa H}{\pa u^i_\al} = \di\frac{\pa L}{\pa u^i_\al} + p_i^\al = 0.$$
This necessary condition (critical point) is also sufficient since the Hamiltonian $H(x,u,p)$ is a concave (and hence incave) function with respect to $u$.
In fact, the Lagrangian $L(x,u)$ is an incave function with respect to $u$, because it has any critical point, and the Hamiltonian $H$
is obtained from $L$ adding a linear term.

We eliminate $p^\al_i$ using the adjoint PDE. It follows the multitime Euler-Lagrange PDEs
$$\di\frac{\pa L}{\pa x^i} -  D_\al\left(\di\frac{\pa L}{\pa x^i_\al}\right)=0.$$
Summarizing, we obtain

{\bf Theorem 6.2}. {\it A $C^2$ solution of the previous optimal control problem is a harmonic map.}

{\bf Remark}. It is possible to extend the notion of harmonic maps to much less
regular maps, which belong to the Sobolev space $H^1(N, M)$ of maps
from $(N, h)$ into $(M, g)$ with finite energy. The above equation is true but only
in the distribution sense and we speak of {\it weakly harmonic maps}.

\section{Minimal volume at constant area as \\optimal control problem}

Suppose that $D$ is a compact set of $R^m = \{(t^1,...,t^m)\}$ with a piecewise smooth $(m-1)$-dimensional boundary $\pa D$.
The volume $\di\int_{D}dt^1...dt^m$ of the domain $D$
is related to the flux of the position vector $t = (t^\al)$
through the closed hypersurface $\pa D$ by the Gauss-Ostragradski formula
$$m\di\int_{D}dt^1...dt^m = \int_{\pa D}\delta_{\al\be} t^\al N^\be \,d\sigma,$$
where $N = (N^\be)$ is the exterior unit normal vector field on $\pa D$. On the other hand,
the area of $\pa D$ is $\di\int_{\pa D}d\sigma$. Introducing a parametrization on $\pa D$, whose domain is $U \subset R^{m-1}$, we have
$d\sigma = ||{\cal N}|| d\eta$, where ${\cal N} = ||{\cal N}|| N$ and $\eta$ is an $(m-1)$-form.

Let us show that {\it of all solids having a given surface area, the
sphere is the one having the greatest volume}. To prove this
statement, we take the {\it normal vector field} ${\cal N}$ as a
{\it control} - a very interesting idea of my PhD student Andreea
Bejenaru - and we formulate the multitime optimal control problem
with isoperimetric constraint
$$\max_{\cal N} \int_{U}\delta_{\al\be} t^\al{\cal N}^\be(t)\, d\eta\,\,\,\,\,
\hbox{subject\,\, to\,\,\,} \di\int_{U}\sqrt{\de_{\al\be} {\cal N}^\al(t){\cal N}^\be(t)}\,\, d\eta = const.$$
Using the Hamiltonian $$H = \delta_{\al\be} t^\al{\cal N}^\be - p\, \sqrt{\de_{\al\be} {\cal N}^\al{\cal N}^\be},\,\,p = const.,$$
the critical point condition, in the multitime maximum principle, gives
$$0= \frac{\pa H}{\pa {\cal N}} = t - p N.$$
Since the Hamiltonian is a concave function of ${\cal N}$, the critical point is a maximum point.
This confirms that $D$ is the sphere $||t||^2 \leq p^2$ in $R^m$.

{\bf Remark.} In the optimal control problems, the Stokes theorem reads: Let $\om$ be a controlled $(p - 1)$ - form, $1\leq p\leq m$,
with compact support on the $p$ - dimensional submanifold $D$. If $\pa D$ denotes the boundary 
of $D$ with its induced orientation, then $\int_{\pa D}\om = \int_{D}d\om,$ where $d$ is the exterior derivative. 
Consequently, we can use as an action functional either 
those defined by the left hand member of the Stokes formula or those defined by the 
right hand member, eventually with new constraints.  

\section{Maximal area constrained by \\a multitime linear flow}

Let $(M,g)$ be a Riemannian manifold and $(TM, G = g +g)$ be its tangent bundle.
Let $\tilde S : \Omega_{0t_0}\to TM$ be an $m$-sheet and suppose that $\tilde S$ is expressed locally by
$x^i = x^i(t), y^i = y^i(t)$ with respect to the induced coordinates $(x^i,y^i)$ in $TM$, and $t \in \Omega_{0t_0}$ as the multitime evolution
parameter. Then the $m$-sheet  $S = \pi \tilde S$ in M is called the projection of the $m$-sheet $\tilde S$ and is represented locally by
$x^i=x^i(t)$. Now let $\xi = \xi^i(x)\frac{\pa}{\pa x^i}$ be a vector field on $M$ tangent to $S$.
Let $\xi^V = \xi^i(x)\frac{\pa}{\pa y^i}$ be the {\it vertical lift} to $TM$. Let $y = y^i\frac{\pa}{\pa y^i}$ be the {\it Liouville vector field} on $TM$.

Let $(\xi(t) = \xi(x(t)),\,\,y(t))$ be the solution of the linear controlled $m$-flow
$$\frac{\pa \xi^i}{\pa t^\al}(t) = A^i_{j\al}(t)\xi^j(t) + A^i_{n+j\al}(t)y^j(t) + B^i_{a}(t) u^a_\al(t)$$
$$\frac{\pa y^i}{\pa t^\al}(t) = A^{n+i}_{j\al}(t)\xi^j(t) + A^{n+i}_{n+j\al}(t)y^j(t) + B^{n+i}_{a}(t) u^a_\al(t)$$
on $TM$. On the tangent bundle we use the area $1$-form
$$\omega = \frac{1}{2}\left(g_{ij}(x) \xi^i(x) \delta y^j (x)- g_{ij}(x) y^i(x)\delta \xi^j(x)\right)$$
$$\delta y^i = dy^i + \Gamma^i_{jk}y^jdx^k, \,\,\delta \xi^i = d\xi^i + \Gamma^i_{jk}\xi^jdx^k,$$
where $\Gamma^i_{jk}$ is the Riemannian connection on $M$.
Giving a closed curve $\tilde C$ on the image of $(\xi(t), y(t))$, we introduce the area
$$\sigma = \frac{1}{2} \int_{\tilde C}\left(g_{ij}(x) \xi^i(x) \delta y^j- g_{ij}(x) y^i\delta \xi^j(x)\right).$$
To find this area, we introduce the pullback
$$\frac{1}{2}\left(g_{ij}(x(t)) \xi^i(t) y^j_\al(t) - g_{ij}(x(t)) y^i(t)x^j_\al(t)\right)dt^\al.$$
Then, the curvilinear integral
$$\sigma (u(\cdot)) =\frac{1}{2} \int_{C}\left(g_{ij}(x(t)) \xi^i(t) y^j_\al(t) - g_{ij}(x(t)) y^i(t)\xi^j_\al(t)\right)dt^\al$$
is the area of a piece from the $m$-surface $(\xi (t),y(t))$ bounded by the curve $C = \pi(\tilde C)$.

We formulate the multitime optimal control problem: $\max_{u(\cdot)} \sigma (u(\cdot))$
{\it subject to the foregoing controlled $m$-flow}. The Hamiltonian $1$-form
$$H_\al = \frac{1}{2}\left(g_{ij}(x(t)) x^i(t) y^j_\al(t) - g_{ij}(x(t)) y^i(t)x^j_\al(t)\right)$$
$$ + p_i(t) \left(A^i_{j\al}(t)x^j(t) + A^i_{n+j\al}(t)y^j(t) + B^i_{a}(t) u^a_\al(t)\right)$$
$$ + q_i(t)\left(A^{n+i}_{j\al}(t)x^j(t) + A^{n+i}_{n+j\al}(t)y^j(t) + B^{n+i}_{a}(t) u^a_\al(t)\right),$$
is linearly affine with respect to the control variables, i.e., this
Hamiltonian $1$-form can be written as
$$H_\al = L_\al + M_a u^a_\al,$$
where
$$M_a(t) = p_i(t)B^i_{a}(t) +q_i(t)B^{n+i}_{a}(t)$$
are the {\it switching functions}. In general, there will be no extremum unless control variables are bounded,
in which case they are expected to be at the boundary of the admissible region (see, linear optimization, simplex method).
Suppose $- 1\leq u^a_\al \leq 1$. When the multitime maximum principle is applied to this type of problems, the optimal control
$u^{*a}_\al$ must satisfies
$$u^{*a}_\al=\left\{\begin{array}{cc}\,1\,\,\,\, \hbox{for}\,\,\, M_a(t) > 0\\ \ \,\hbox{?\,\, for}\,\, M_a (t)= 0 \\ \ - 1\,\,\,\, \hbox{for}\,\, M_a(t) < 0,\end{array}\right.$$
for each $\al = 1,...,m$. This optimal control is discontinuous since each component
jumps from a minimum to a maximum and vice versa in response to each change
in the sign of each $M_a(t)$ (switching functions). The optimal control $u^{*a}_\al$ is called a {\it bang bang control}.

\section{Conclusions}

This paper refers to basic problems in control the partial differential equations of differential geometry.
The original results are meaningful and useful for explaining many real world phenomena
based on optimal controlled multitime evolutions. They shows that some dreams
issuing from the papers of Lev S. Pontryaguin, Lawrence C. Evans and Jacques-Louis Lions are now partially covered
by the multitime maximum principle. Of course, to pass from the previous local theory to a global one,
we need another more flexible formulation of a smooth multitime optimal control problem involving two
Riemannian manifolds.

{\bf Acknowledgements:} Partially supported by University Politehni-ca of Bucharest,
and by Academy of Romanian Scientists, Bucharest, Romania.

Some ideas were presented at The International Conference of Differential Geometry and Dynamical Systems
(DGDS-2010), 25-28 August 2010, University Politehnica of Bucharest, Bucharest, Romania.

Dedicated to Lawrence C. Evans, Lev S. Pontryaguin and Jacques-Louis Lions for their seminal contributions to optimal control theory.

Constantin Udri\c ste, University POLITEHNICA of Bucharest, Faculty of Applied Sciences,
Department of Mathematics-Informatics I, Splaiul Independentei 313, Bucharest 060042, Romania \\udriste@mathem.pub.ro;
\\anet.udri@yahoo.com;

\end{document}